\documentclass[12pt,a4paper]{article}

\usepackage{amsmath,amssymb}
\usepackage{amsthm}
\usepackage[english]{babel}
\usepackage{multicol}
\usepackage{hyperref}
\usepackage{url}
\usepackage{graphicx}
\usepackage{latexsym}
\usepackage{tikz}
\tikzset{
	every node/.style={circle, inner sep=2pt}
}
\usetikzlibrary{matrix,arrows,calc}

\usepackage{amsmath,amsfonts,amssymb,latexsym,graphics,epsfig,url}
\usepackage{color}
\usepackage{amsthm}
\usepackage[english]{babel}
\usepackage{graphicx}
\usepackage{soul}
\usepackage[utf8]{inputenc}

\newtheorem{theorem}{Theorem}[section]
\newtheorem{proposition}[theorem]{Proposition}
\newtheorem{lemma}[theorem]{Lemma}
\newtheorem{corollary}[theorem]{Corollary}

\newtheorem{problem}[theorem]{Problem}

\newtheorem{definition}[theorem]{Definition}

\input amssym.def
\newsymbol\rtimes 226F
\newfont{\nset}{msbm10}

\def\A{{\mbox {\boldmath $A$}}}
\def\B{{\mbox {\boldmath $B$}}}

\def\D{{\mbox {\boldmath $D$}}}

\def\G{\Gamma}
\def\I{{\mbox {\boldmath $I$}}}
\def\J{{\mbox {\boldmath $J$}}}

\def\L{{\mbox {\boldmath $L$}}}
\def\M{{\mbox {\boldmath $M$}}}

\def\P{{\mbox {\boldmath $P$}}}

\def\A{{\mbox {\boldmath $A$}}}

\def\matrix0{{\mbox {\boldmath $O$}}}

\def\S{\mbox{\boldmath $S$}}

\def\vv{{\mbox{\boldmath $v$}}}

\def\vec0{\mbox{\bf 0}}

\def\diag{\mathop{\rm diag }\nolimits}

\def\mod{\mathop{\rm mod }\nolimits}

\def\Ker{\mathop{\rm Ker }\nolimits}
\def\rank{\mathop{\rm rank }\nolimits}
\def\tr{\mathop{\rm tr }\nolimits}

\def\spec{\mathop{\rm sp }\nolimits}

\def\G{\Gamma}

\begin{document}
	
	\title{On token signed graphs
		\thanks{The research of C. Dalf\'o and M. A. Fiol has been supported by
			AGAUR from the Catalan Government under project 2021SGR00434 and MICINN from the Spanish Government under project PID2020-115442RB-I00.
			The research of M. A. Fiol was also supported by a grant from the  Universitat Polit\`ecnica de Catalunya with references AGRUPS-2022 and AGRUPS-2023.}
	}
	\author{C. Dalf\'o$^a$, M. A. Fiol$^b$, and E. Steffen$^c$\\
		\\
		{\small $^a$Departament. de Matem\`atica}\\ {\small  Universitat de Lleida, Lleida/Igualada, Catalonia}\\
		{\small {\tt cristina.dalfo@udl.cat}}\\
		{\small $^{b}$Departament de Matem\`atiques}\\ 
		{\small Universitat Polit\`ecnica de Catalunya, Barcelona, Catalonia} \\
		{\small Barcelona Graduate School of Mathematics} \\
		{\small  Institut de Matem\`atiques de la UPC-BarcelonaTech (IMTech)}\\
		{\small {\tt miguel.angel.fiol@upc.edu}}\\
		{\small $^{c}$Department of Mathematics}\\
		{\small Paderborn University, Warburger Str. 100,}\\
		{\small 33098 Paderborn, Germany }\\
		{\small {\tt es@uni-paderborn.de} }\\
	}

	\maketitle

	\begin{abstract}
		We introduce the concept of a $k$-token  signed graph and study some of its combinatorial and algebraic properties. We prove that two switching isomorphic signed graphs have switching isomorphic token graphs.
		Moreover, we show that the Laplacian spectrum of a balanced signed graph is contained in the Laplacian spectra of its  $k$-token signed graph. Besides, we introduce 
		and study the
		\emph{unbalance level} of a signed graph, which is a new parameter that measures how far a signed graph is from being balanced. Moreover, we 
		study the relation between the frustration index and the
		unbalance level of signed graphs and their token signed graphs.
	\end{abstract}
	
	\textbf{Keywords:} Token graph, Signed graph, Laplacian spectrum, Adjacency spectrum, Frustration index.
	
	%%%%%%%%%%%%%%%%%%%%%%%%%%%%%%%%%%%%%%%%%%%%%%%%%%%%%%%%%%%%%%%%%%%	
	\section{Introduction}
	% Given a graph $G=(V,E)$ on $n$ vertices, its $k$-token graph $F_k(G)$ has a vertex set consisting of all the ${n\choose k}$ $k$ indistinguishable tokens placed in different vertices of $G$, and two vertices are adjacent if one is obtained from the other by moving one token along an edge of $G$.
	Given a graph $G=(V,E)$ on $n$ vertices, its $k$-token graph $F_k(G)$ has ${n\choose k}$ vertices corresponding to configurations of $k$ indistinguishable tokens placed at distinct vertices of $G$, and two
	configurations are adjacent whenever one configuration can be reached from the other
	by moving one token along an edge from its current position to an unoccupied vertex.
	Alternatively, we can view the vertices of $F_k(G)$ as $k$-subsets of $V$, two vertices $A$ and $B$ being adjacent if their symmetric difference $A\Delta B$ consists of two adjacent vertices in $G$. Token graphs were also called {\em symmetric $k$-th power of a graph} by Audenaert, Godsil, Royle, and Rudolph \cite{agrr07}, and {\em $k$-tuple vertex graphs} in  Alavi, Behzad, Erd\H{o}s, and Lick \cite{abel91}.
	Token graphs have applications in physics (in quantum mechanics) and in the graph isomorphism problem (because, usually, invariants of the $k$-token $F_k(G)$ are also invariants of $G$). For more details, see again \cite{agrr07}.
	Some properties of token graphs were studied by Fabila-Monroy, Flores-Pe\~{n}aloza, Huemer,  Hurtado,  Urrutia, and  Wood \cite{ffhhuw12}.
	
	In a signed graph, every edge is given a positive or negative sign. Thus,
	signed graphs are a powerful tool for modeling and analyzing various scenarios where relationships can have both positive and negative aspects. Applications of signed graphs include Social Network Analysis, Psychology and Sociology, Natural Sciences,
	Game Theory, Complex Systems, etc. For the notation and main properties on signed graphs, see, for example, Belardo, Cioab\v{a}, Koolen, and Wang \cite{bckw18}.
	
	In this paper, we `merge' both token and signed graph concepts by defining the $k$-token signed graphs and studying some of their combinatorial and algebraic properties. In particular, we prove that two switching equivalent signed graphs have switching equivalent token graphs.
	Moreover, we show that the Laplacian spectrum of a balanced signed graph is contained in the Laplacian spectra of its  $k$-token signed graph.
	
	This paper is structured as follows. In the next section, we give the basic concepts and results. In Section \ref{sec:token}, we generalize the concept of token graphs to token signed graphs and present some properties of the latter. In section \ref{sec:frust-unbal}, we discuss two measures of `unbalance' in a signed graph. The first one is the well-known frustration index, and the second one is a new spectral measure called
	the unbalance level. Some examples show that this new measure has a very
	good discernment capacity between different switching isomorphic classes of graphs. In Section \ref{sec:sign-symm}, we show that $k$-tokens graphs preserve sign-symmetry. In Section \ref{sec:Lapl}, we define the signed $(n,k)$-binomial matrix $\B$. This matrix allows us to prove that, given two different token graphs $F_{k_1}$ and $F_{k_2}$, $\B\L_{k_1}=\L_{k_2}\B$. Finally, in Section \ref{sec:compl}, we define the complement of a balanced signed graph. Then, given the Laplacian matrices of the $k$-token graphs of a balanced signed graph $\L_k$ and its signed complement graph 
	$\overline{\L}_k$, we show that $\L_k\overline{\L}_k=\overline{\L}_k \L_k$.

	%%%%%%%%%%%%%%%%%%%%%%%%%%%%%%%%%%%%%%%%%%%%%%%%%%%%%%%%%
	\section{Basic concepts and results}
	
	Let us first recall some definitions and basic results about signed graphs. For more details, see, for instance, Zaslavsky \cite{z82,z10}.
	
	Let $G$ be a graph with vertex set $V=V(G)$ and edge set $E=E(G)$.
	Then, the {\em signed graph} $\G=(G,\sigma)$ is a graph $G$ together with a function $\sigma: E\rightarrow \{+1,-1\}$ called the {\em signature} of $G$. If $\sigma(e) = 1$ (respectively, $\sigma(e) = -1$) for every edge $e$, then
	$\sigma$ is called the {\em all-positive} (respectively, {\em all-negative}) signature and 
	$(G,\sigma)$ is called an {\em all-positive} (respectively, {\em all-negative}) signed graph. 
	We consider simple graphs and an edge $e$, which is incident to vertices $u$ and $v$, is
	also denoted by $uv$. 
	We represent positive and negative edges in the figures by thin and thick lines, respectively (see Figure \ref{fig1}).
	
	Given a signed graph $\G=(G,\sigma)$ and a vertex subset $U$, the {\em switching at $U$}  yields the graph $\G^U=(G,\sigma^U)$ obtained by reversing the sign of all the edges of $\partial_{\G}(U)$, which denotes the set of edges between $U$ and $V(G)\setminus U$. We say that $\G$ and $\G'=\G^U$ are {\em (switching) equivalent}, denoted by $\G\sim \G'$, or also that the two signatures $\sigma$ and
	$\sigma^U$ are equivalent. See an example in Figure \ref{fig2} (left).
	
	The {\em sign of a set} $E$ of edges in $\G$ is the product of the signs of its elements,
	and it is called {\em positive} (respectively, {\em negative}) if its sign is positive (respectively, negative). Clearly, switching does not change the sign of a cycle. 
	A signed graph is {\em balanced} if all its cycles are positive, and {\em unbalanced} otherwise.
	Harary \cite{Harary53} characterized balanced signed graphs.
	
	\begin{theorem}[\cite{Harary53}] \label{Theo:Harary}
		A signed graph $(G,\sigma)$
		is \emph{balanced} if and only if $V(G)$ can be partitioned into two
		sets, one of which might be empty, such that every positive edge connects two vertices of the same set and 
		every negative edge connects two vertices of different sets. 
	\end{theorem}
	
	The following corollary is equivalent to Theorem \ref{Theo:Harary}.
	
	\begin{corollary} \label{Cor:charact_balanced}
		A signed graph $(G,\sigma)$ is balanced if and only if it is switching equivalent
		to the all-positive signed graph $+G$. 
	\end{corollary}
	
	The {\em (signed) adjacency matrix} $\A=(a_{uv})$ of a signed graph $\G=(G,\sigma)$ 
	has entries $a_{uv}=\sigma(uv)$ if $uv$ is an edge of $G$, and $a_{uv}=0$, otherwise. The Laplacian matrix of $\G$ is $\L=\D-\A$, where $\D$ is the diagonal matrix of the degrees of $G$.
	Switching equivalent signed graphs, $\G=(G,\sigma)$ and $\G'=(G,\sigma')$, have similar adjacency matrices and, hence, they are cospectral.
	Indeed, two signatures $\sigma$ and $\sigma'$ are equivalent if there exists an  $n\times n$ diagonal matrix $\S$,  with entries $s_u=+1$ if $u\in U$, and $s_u=-1$ otherwise, such that, if $\A$ and $\L$ [respectively, $\A'$ and $\L'$] are the adjacency and Laplacian matrices of $\G$ [respectively, $\G'$], we have
	\begin{equation}
		\A'=\S\A\S\quad\mbox{ and }\quad \L'=\S\L\S.
	\end{equation}
	
	Two simple signed graphs $\G= (G,\sigma)$ and 
	$\G' = (H,\tau)$ are {\em isomorphic}, denoted by $\G\cong \G'$, if there is a 
	bijection $\phi: V(G) \to V(H)$ such that  
	$e=uv\in E(G)$ if and only if $\phi(u)\phi(v) \in E(H)$ and  $\tau(\phi(u)\phi(v)) = \sigma(e)$. The mapping $\phi$ is called a {\em (sign-preserving) isomorphism} between 
	$(G,\sigma) $ and $(H,\tau)$. In particular, an {\em automorphism} of 
	a signed graph is an isomorphism of the graph to itself. In terms of their adjacency or Laplacian matrices,  the signed graphs $\G$ and $\G'$ are isomorphic when there exists a permutation matrix $\P$ such that
	\begin{equation}
		\A'=\P\A\P^{\top}\quad\mbox{ and }\quad \L'=\P\L\P^{\top}.
	\end{equation}
	
	Two signed graphs $\G=(G,\sigma)$ and $\G'=(G',\sigma')$ are {\em switching isomorphic,} denoted by $\G\simeq \G'$,  
	if $\G'$ is isomorphic to $\G^U$ for a subset $U$ of $V(G)$. 
	Note that, if $\G \simeq \G'$ and $\phi$ is an isomorphism
	between $\G^U$ and $\G'$, then $\phi^{-1}$ is an isomorphism between $\G$ and $(\G')^{\phi[U]}$.
	Thus, we can say that $\G$ and $\G'$ are switching isomorphic if one of them is isomorphic to
	a switching equivalent signed graph of the other (the symmetry property of the relation $\simeq$).
	In terms of the adjacency and Laplacian matrices of $\G$ and $\G'$, this means that there exist a diagonal matrix $\S$, and a permutation matrix $\P$, such that
	\begin{equation}
		\A'=\P\S\A\S\P^{\top}\quad\mbox{ and }\quad \L'=\P\S\L\S\P^{\top}.
		\label{eq:switch-isomof}
	\end{equation}
	Then, since $\P\P^{\top}=\P^{\top}\P=\S^2=\I$, from \eqref{eq:switch-isomof}, we get $\S\P^{\top}\A'\P\S=\A$ or
	$$
	\A=\P^{\top}(\P\S\P^{\top})\A'(\P\S\P^{\top})\P=\P^{\top}\S'\A'\S'\P,
	$$
	where $\S'=\P\S\P^{\top}$ and $\P^{\top}(\cdots)\P$ correspond, respectively, to the above $\phi(U)$ and $\phi^{-1}$. 
	
	Thus, switching isomorphic signed graphs are also cospectral.
	It is easy to see that a signed graph on $n$ vertices, $m$ edges, and $c$ components 
	has $2^{m-n+c}$ pairwise different signatures with respect to switching equivalence. Since any two switching equivalent signed graphs are also 
	switching isomorphic, there are at most that many classes of switching isomorphic signatures on a graph $G$. However, it can be much less. For instance, the 
	Petersen graph has $2^6$ pairwise different signatures with respect to switching equivalence but, as shown by Zaslvsky \cite[Th. 5.1]{z12}, only $6$ pairwise different signatures with respect to switching isomorphism, see Table \ref{tab:classesP}.
	
	The {\em negation} $-\G$ of a signed graph $\G=(G,\sigma)$ is the graph obtained by reversing all the signs of $\G$. That is, $-\G=(G,-\sigma)$.
	A signed graph $\G$ is {\em sign-symmetric} if it is switching isomorphic to its negation. 
	%\blue{$-\G$. Then, a sign-symmetric graph must have 
		
		Now, let us turn our attention to token graphs, already defined in the Introduction.
		The following concepts and results about the Laplacian matrices of token graphs were given by Dalf\'o, Duque, Fabila-Monroy, Fiol, Huemer, Trujillo-Negrete, and Zaragoza Mart\'inez in \cite{ddffhtz21}.
		For some integers $k_1$ and $k_2$, with $1\le k_1<k_2<n$, the {\em $(n;k_1,k_2)$-binomial matrix} $\B$ is a ${n\choose k_2}\times{n\choose k_1}$ matrix whose rows and columns are indexed by the $k_2$-subsets $A\subset [n]=\{1,\ldots,n\}$ and $k_1$-subsets $X\subset [n]$, respectively, with entries
		$$
		(\B)_{AX}=\left\{
		\begin{array}{ll}
			1 & \mbox{if $X\subset A$,}\\
			0 & \mbox{otherwise}.
		\end{array}
		\right.
		$$
		\begin{theorem}[\cite{ddffhtz21}]
			\label{theo:BL1=LB2}
			Let $G$ be a graph on $n$ vertices,
			with $k_1$-token and $k_2$-token graphs $F_{k_1}(G)$ and $F_{k_2}(G)$, where $1\le k_1\le k_2<n$. Let $\L_{k_1}$ and $\L_{k_2}$ be the respective Laplacian matrices, and $\B$ the $(n;k_1,k_2)$-binomial matrix. Then,
			$$
			\B\L_{k_1}=\L_{k_2}\B.
			$$  
		\end{theorem}
		\begin{theorem}[\cite{ddffhtz21}]
			\label{theo:LnoL=noLL}
			Given a graph $G$ and its complement $\overline{G}$, the Laplacian matrices $\L=\L(F_k(G))$, and $\overline{\L}=\L(F_k(\overline{G}))$ of their $k$-token graphs commute:
			$$
			\L\overline{\L}=\overline{\L}\L.
			$$  
		\end{theorem}
		
		%%%%%%%%%%%%%%%%%%%%%%%%%%%%%%%%%%%%%%%%%%%%%%%%%%%%%%%%%%%%%%
		\section{Token signed graphs}
		\label{sec:token}
		
		The concept of token graphs generalizes in a natural way to token signed graphs, as shown in the following definition.
		
		\begin{definition}
			Let $\Gamma = (G,\sigma)$ be a simple signed graph and 
			$k \geq 1$ be an integer. Let $(F_k(\Gamma),\sigma_k)$ be the signed graph with vertex set $V(G) \choose k$, where two vertices $A$ and $B$ of $F_k(\Gamma)$ are adjacent if $A \Delta B = \{a,b\}$ with $a\in A$, $b \in B$, $ab \in E(G)$, and $\sigma_k(A,B) = \sigma(ab)$.    
		\end{definition}
		
		Equivalently, the vertices of $F_k(\G)$ correspond to $k$ indistinguishable tokens placed in different vertices of $G$,
		and there is an edge between $A$ and $B$, with sign $\sigma_k(A,B)=s$, if $B$ is obtained from $A$ by moving one token from $A$ to $B$ along an edge $e\in E(G)$  with sign $\sigma(e)=s$.
		
		Let $\G = (G,\sigma)$ be a signed graph.
		The signature $\sigma_k$ of $(F_k(\G),\sigma_k)$ is naturally determined by $\sigma$. For this reason, and also to avoid overloading of notation, we denote the signed
		token graph $(F_k(\G),\sigma_k)$ by $F_k(\G)$
		in the following. 
		
\begin{figure}[t]
			\begin{center}
				\includegraphics[width=10cm]{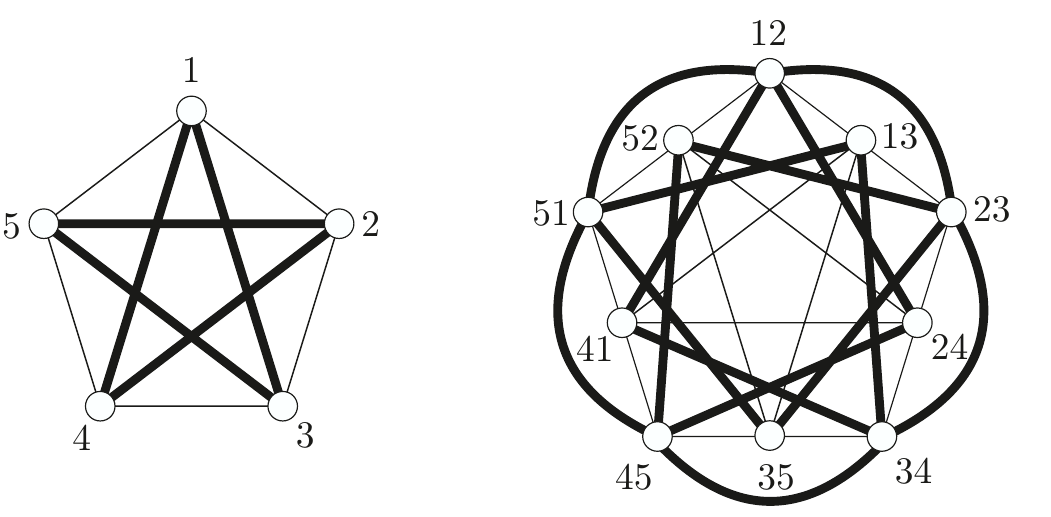}
				%\vskip -1.5cm
				\caption{A signed $K_5$ and its 2-token signed graph.} 
				%, \textcolor{blue}{which is balanced if $V_1=\{12,23,34,45,51\}$ and $V_2=\{13,24,35,42,52\}$}. 
				\label{fig1}
			\end{center}
		\end{figure}

\begin{figure}[!ht]
				\begin{center}
					\includegraphics[width=10cm]{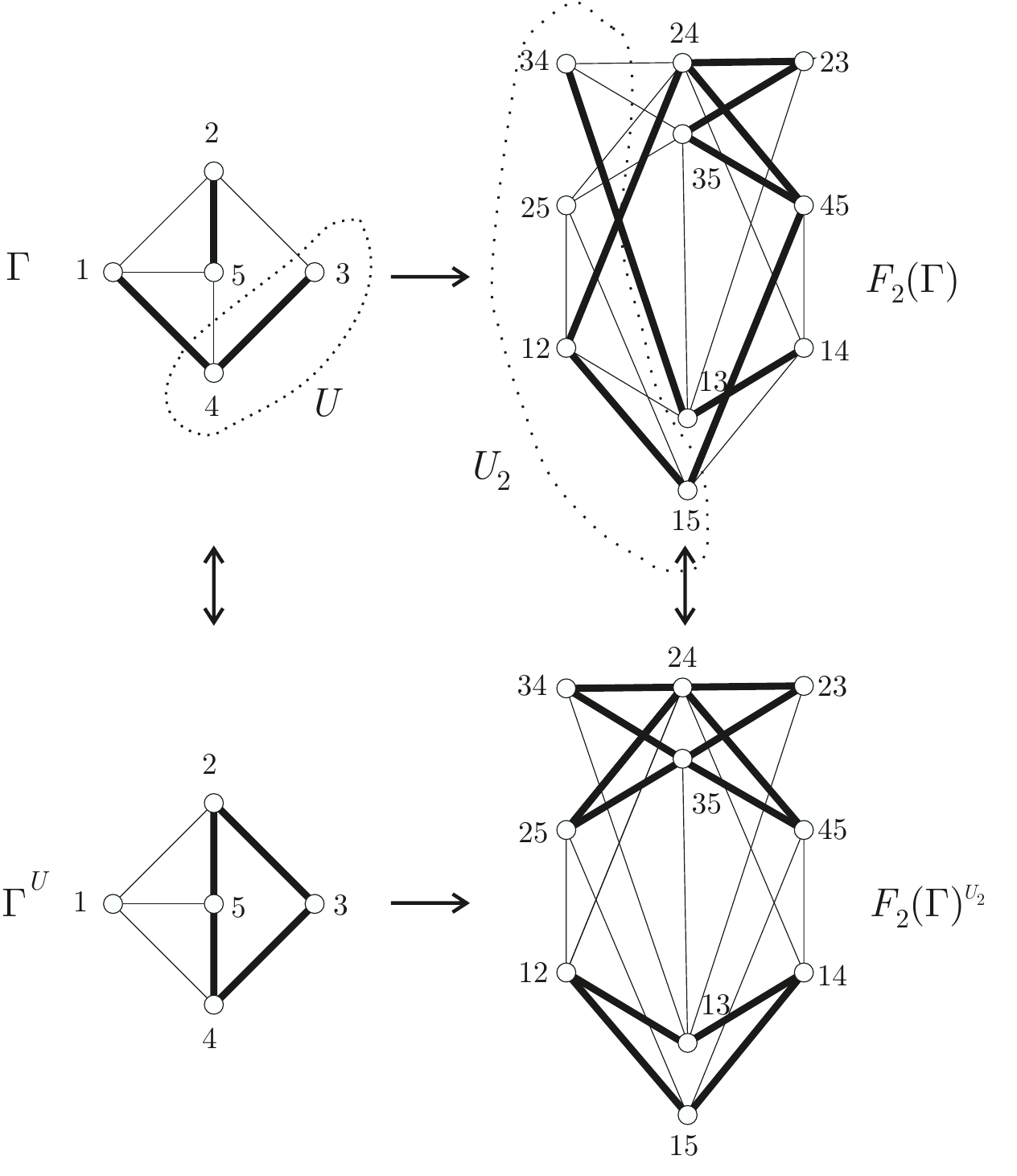}
					%\vskip -1.5cm
					\caption{Token graphs preserve switching equivalence.}
					\label{fig2}
				\end{center}
			\end{figure}

		In particular, notice that
		$F_1(\G)=\G$ and, by symmetry,
		$F_k(\G)=F_{n-k}(\G)$. Moreover,
		$F_k(K_n)$ is the signed Johnson graph $(J(n,k),\sigma)$.
		For example, Figure \ref{fig1} shows a signed Johnson graph $J(5,2)$ as a 2-token graph of a signed complete graph $K_5$.
		
		The following result shows how positive and negative signs in a signed graph $\G$ translate into corresponding signs in $F_k(\G)$.
		\begin{lemma}
			\label{lem:cycles}
			Let $\G=(G,\sigma)$  be a signed graph with $n$ vertices and $m$ edges. Let $F_k(\G)$ be its $k$-token graph for some $k\le n/2$. 
			\begin{itemize}
				\item[$(i)$] 
				If $\G$ has $m^+$ positive edges and $m^-$ negative edges, then $F_k(\G)$ has ${n-2\choose k-1}m^+$ positive edges and ${n-2\choose k-1}m^-$ negative edges.
				\item[$(ii)$] 
				If $\G$ has a positive (respectively, negative) cycle of length $p$, for $3\le p\le n$,
				and $k'\le k$ is an integer satisfying $k+p-n\le k'<p$, then $F_k(\G)$ has ${n-p\choose k-k'}$ positive (respectively, negative) cycles of length $p$.
			\end{itemize}
		\end{lemma}
		
		\begin{proof}
			The statement $(i)$ is clear when the adjacencies of $F_k(\G)$ are defined in terms of token movements.
			Then, let us prove $(ii)$ by using the same approach. Let $0,1,\ldots,p-1$ denote the vertices of the cycle. First note that, if $k'\le \min\{k,p\}$
			and $k-k'\le n-p$, that is, $k'\ge p+k-n$, we can place $k'$ tokens on the cycle (see Table \ref{tab:lemma}, on the left, with $n=7$ and $k=3$) with a given order 
			$(\mod p)$, 
			say, $0\le i_1<i_2<\cdots <i_{k'}\le p-1$. This gives ${n-p\choose k-k'}$ possible places for the remaining tokens. Then, we only need to prove that the given $k'$ tokens on the cycle correspond to a vertex of a positive or negative $p$-cycle in $F_k(\G)$. With self-explanatory notation, the movement of the tokens (in a positive circular sense modulo $p$) are as shown on the right of Table \ref{tab:lemma}. 
			Notice that, we have done $p$ movements of tokens along all the $p$ edges of the cycle. Thus, this corresponds to a $p$-cycle in $F_k(\G)$ of the same sign as the $p$-cycle in $\G$.
		\end{proof}
		
		\begin{table}
			\begin{center}
				\begin{tabular}{|c c c|}
					\hline
					& $n=7$, $k=3$ &  \\
					\hline
					$p$ & $p+k-n$ & $k'$\\
					\hline
					$3$ & $-1$ & $1,2$\\
					$4$ & $0$ & $1,2,3$\\
					$5$ &  $1$ & $1,2,3$\\
					$6$ &  $2$ & $2,3$\\
					$7$ & $3$ & $3$\\
					\hline
				\end{tabular}\qquad\quad
				\begin{tabular}{|lcl|}
					\hline
					$i_{k'}$ & $\rightarrow$ & $i'_{k'}=i_1-1$\\[0.05cm]
					$i_{k'-1}$ & $\rightarrow$ & $i'_{k'-1}=i_{k'}$\\[0.05cm]
					$i_{k'-2}$ & $\rightarrow$ & $i'_{k'-2}=i_{k'-1}$\\[0.05cm]
					& \vdots &  \\[0.05cm]
					$i_1$ &  $\rightarrow$ & $i'_1=i_2$\\[0.05cm]
					$i'_{k'}$ & $\rightarrow$ & $i''_{k'}=i_1$\\[0.05cm]
					\hline
				\end{tabular}
				\caption{Scheme of the proof of Lemma \ref{lem:cycles}}
				\label{tab:lemma}
			\end{center}
		\end{table}
		
		In the following result, we show that the operation of taking $k$--tokens preserves the switching equivalence and switching isomorphism.
		
		%	\begin{theorem} \label{Theo: switch_equivalent G --> F_k(G)}
			%		\label{theo:switch-equiv}
			%		 For every signed graph $\G = (G,\sigma)$ and
			%		every $U \subseteq V(G)$,
			%	there exist a subset $U_k \subseteq V(F_k(\G))$ such that
			%	$$
			%	F_k(\G^U)=F_k(\G)^{U_k}.
			%	$$
			%\end{theorem}
			
			\begin{theorem} 
				%\label{Theo:switch_equivalentG-->F_k(G)}
				\label{theo:switch-equiv}
				Let $1 \leq k \leq n$ be integers and $\G$ and $\G'$ be signed graphs 
				of order $n$. 
				\begin{enumerate}
					\item[$(i)$] 
					If $\G$ and $\G'$ are switching equivalent, then $F_k(\G)$ and $F_k(\G')$ are switching equivalent.  
					\item[$(ii)$] 
					If $\G$ and $\G'$ are switching isomorphic, then $F_k(\G)$ and $F_k(\G')$ are switching isomorphic.
				\end{enumerate}
			\end{theorem}
			
			\begin{proof}
				Let $\G = (G,\sigma)$, $\G' = (G,\sigma')$ and 
				$V(G) = \{v_1, \dots, v_n\}$.
				
				$(i)$ We have to prove that for every $U \subseteq V(G)$,
				there exist a subset $U_k \subseteq V(F_k(\G))$ such that
				$$
				F_k(\G^U)=F_k(\G)^{U_k}.
				$$
				
				Define $U_k = \{A \colon A \in V(F_k(\G)) \text{ and } \vert A \cap U \vert \text{ is even}\}$.
				
				Let $\sigma'$ be the signature of $\G^U=\G'$, 
				$\sigma_k$ be the signature of $F_k(\G)$,
				$\sigma_k'$ be the signature of $F_k(\G^U)$, and let
				$\sigma_k^*$ be the signature of $F_k(\G)^{U_k}$. 
				Let $xy \in E(G)$ and $AB$ be an edge of $F_k(\G)$ with $A \Delta B = \{x,y\}$.
				Note that $\sigma_k(AB) = \sigma(xy)$ and $\sigma_k'(AB)=\sigma'(xy)$.
				We have to show that $\sigma_k^*(AB) = \sigma_k'(AB)$. 
				
				If $x,y$ are both in $U$ or both in $V(G)\setminus U$, then  
				$\sigma(xy) = \sigma'(xy)$ and consequently, $\sigma_k(AB)=\sigma_k'(AB)$. 
				Furthermore, either both, $A$ and $B$, are in $U_k$ or
				both are in $V(F_k(\G))\setminus U_k$ and, therefore,  $\sigma_k^*(AB) = \sigma_k(AB) = \sigma_k'(AB)$.
				
				Thus, we may assume that $x \in U$ and $y \not \in U$. Then, precisely one of $A$
				and $B$ is an element of $U_k$, say $A$. Then 
				$\sigma_k^*(AB) = - \sigma_k(AB) = -\sigma(xy)= \sigma'(xy) = \sigma_k'(AB)$. 
				
				$(ii)$
				Since $\G$ and $\G'$ are switching isomorphic, 
				there is $U \subseteq V(G)$
				such that $\G'$ is isomorphic to $\G^U$. Let $\phi \colon V(\G') \rightarrow V(\G')$ be the corresponding automorphism. 
				So, $\G$ and $\G^U$ are switching equivalent and, thus, by $(i)$, 
				there is $U_k \subseteq V(F_k(\G))$ such that $F_k(\G^U)=F_k(\G)^{U_k}$.
				We have to show that $F_k(\G')$ is isomorphic to $F_k(\G)^{U_k}$.
				
				Let $A = \{a_1, \dots, a_k\} \in V(F_k(\G'))$. Then, $\phi' \colon V(F_k(\G') \rightarrow V(F_k(\G')$ with 
				$\phi'(A) = \{\phi(a_1), \dots, \phi(a_k)\}$ is an automorphism on
				$F_k(\G')$ that maps $F_k(\G')$ to $F_k(\G)^{U_k}$. Thus, $F_k(\G)$ and $F_k(\G')$ are switching isomorphic.
			\end{proof}
			
			%%%%%%%%%%%%%%%%%%%%%%%%%%%%%%%%%%%%%%%%%%%%%%%%%%%%%%%%%%%%
			\section{Frustration index and unbalance level}
			\label{sec:frust-unbal}
			
			In this section, we discuss two measures of `unbalance' in a signed graph. The first one is the well-known frustration index, which has received a lot of attention in the literature. The second one is a new spectral measure called the \textit{unbalance level}. Some examples show that this new measure has a very good discernment capacity between different switching isomorphic classes of graphs.
			
			The {\em frustration index} of a signed graph $\G=(G,\sigma)$, denoted by $l(\G)$ or $l(G,\sigma)$, is the minimum number
			of negative edges that we need to remove to obtain a balanced graph. 
			%A signed graph $(G,\sigma)$ %with $l(G,\sigma) = k$ is $called $k$-frustrated 
			Thus, the frustration index gives a measure of how far a signed graph is from being balanced.
			
			Concerning the frustration index of token graphs, we have the following result. 
			
			\begin{proposition}
				Let $\G=(G,\sigma)$ be a signed graph of order $n$ and frustration index $l(\G)$. Then, the frustration index of its $k$-token graph satisfies
				\begin{equation}
					l(\G) \le l(F_k(\G))\le {n-2\choose k-1} l(\G).
					\label{eq:bounds-l}
				\end{equation}
			\end{proposition}
			\begin{proof} The frustration index does not change under switching, and it is well-known that 
				there is $U \subseteq V(G)$ such that $\G^U$ has precisely
				$l(\G)$ negative edges; see, for instance, Lemma 1.1 in Cappello and Steffen \cite{cs22}.
				By Lemma \ref{lem:cycles}$(i)$, every negative edge of $\G$ gives rise to
				${n-2\choose k-1}$ edges in $F_k(\G)$ and, thus, 
				by Theorem \ref{theo:switch-equiv}, 
				we get the upper bound in \eqref{eq:bounds-l}.\\
				To prove the lower bound, consider the set ${\cal S}(e)$ of ${n-2\choose k-1}$ edges induced by every edge $e\in E(G)$. Now, reasoning by contradiction,  assume that $l(F_k(\G))=t<l(G)$. As before, there exists a vertex set $U_k\subset V(F_k(G))$ such that $F_k(\G)^{U_k}$ has exactly $t$ negative edges ${\cal E}_i\in {\cal S}(e_i)$ for some  $e_i\in E(G)$, with $i=1,\ldots,t$. Hence, $F_k'=F_k(\G)^{U_k}\setminus \{{\cal E}_1,\ldots,{\cal E}_t\}$ would be balanced. Moreover, the multiset $\{e_1,\ldots,e_t\}$ has at most $t$ different elements, and  $\G\setminus \{e_1,\ldots,e_t\}$ would be also  a balanced graph. Otherwise, any negative cycle would induce some negative cycle in $F_k'$ by Lemma \ref{lem:cycles}$(ii)$. Consequently, $l(\G)\le t$, a contradiction.
			\end{proof}
			As a consequence of Theorem \ref{theo:switch-equiv}
			we obtain the following result. 
			
			\begin{corollary} \label{coro:F_kbalanced}
				Let $\G=(G,\sigma)$ be a balanced signed graph on $n$
				vertices. Then, for every integer $k$, with $1\le k\le n-1$, the  $k$-token signed graph $F_k(\G)$ is also balanced.
			\end{corollary}
			\begin{proof}
				By Corollary \ref{Cor:charact_balanced}, $\G$ is switching equivalent to the all-positive signed graph $\G'$. Thus, $F_k(\G')$ is all positive and switching equivalent
				to $F_k(\G)$ by Theorem \ref{theo:switch-equiv}. Hence,
				$F_k(\G)$ is balanced. 
			\end{proof}
			% \blue{
				% \subsection{Unbalance level}
				% }
			In the second part of this section, we introduce a new measure of balance, which is based on the spectra of the signed graphs considered. The idea is to give a number related to the fraction of negative cycles with respect to the total. In terms of the traces of the matrices $\A^+$ and  $\A$, this can be done by using the parameter
			\begin{equation}
				\ell_m(\G)=\frac{\sum_{r=0}^m\left[\tr(\A^+)^r-\tr \A^r\right]}
				{\sum_{r=0}^m\left[\tr(\A^+)^r+|\tr \A^r|\right]}
				\label{ell-m}
			\end{equation}
			for some integer $m$ related to the number $n$ of vertices of the signed graph. Then, to have a record of both the even and odd cycles, computer evidence leads us to the following definition.
			
			\begin{definition}
				Let $\G=(G,\sigma)$ be a signed graph on $n$ vertices. Let $\A$ and $\A^+$ be the adjacency matrices of $\G$ (with signed matrix) and $G$ (with unsigned matrix), respectively.
				The {\em spectral unbalance level}
				of $\G$ is a number in the interval $[0,1]$ defined in terms of \eqref{ell-m} as
				$$
				\ell(\G)=\max\{\ell_{n-1}(\G),\ell_n(\G)\}.
				$$
			\end{definition}
			As commented, this is a spectral measure since the traces of $\A^r$ and $(\A^+)^r$ in \eqref{ell-m} can be computed by using their respective eigenvalues. Indeed, if $\A$ has eigenvalues $\lambda_1,\ldots,\lambda_n$, then 
			$$
			\tr \A^r=\sum_{i=0}^n \lambda_i^r,
			$$
			and, similarly, for $\tr(\A^+)^r$.
			Moreover, since a balanced signed graph $\G=(G,\sigma)$ is switching equivalent to (the all-positive) graph $+G$, we have that $\ell(\G)=0$ if and only if $\G$ is balanced.
			In the other extreme, the highest `degree of unbalance' would be when $\ell(\G)=1$.
			In Table \ref{tab:all(-Cn)}, we show the unbalanced level of the cycles $C_n^-$ with only one negative edge and that of the all-negative cycles $-C_n$. Looking at the table, the following remarks are worth mentioning:
			\begin{enumerate}
				\item 
				As $n$ increases, the unbalanced level of both $C_n^-$ and $-C_n$ tends to zero. 
				\item
				For $n$ even, $\ell(-C_n)=0$ since $-C_n$ is switching equivalent to the all-positive cycle $C_n$.
				\item 
				For $n$ odd, $\ell(C_n^-)=\ell(-C_n)$ since, in this case, $C_n^-$ and $-C_n$ are switching equivalent.
			\end{enumerate}
			
			Following the examples, Table \ref{tab:all(-Kn)} shows the unbalanced level of the complete graphs $K_n^-$ with only one negative edge and that of the all-negative complete graphs $-K_n$.
			As in the case of cycles, we can easily compute such unbalanced levels by using the spectra of $K_n$ (that is, $\{n,-1^{n-1}\}$), $-K_n$ (that is, $\{-n,1^{n-1}\}$), and $K_n^-$ (see Table \ref{tab:sp(Kn^-)}).
			Here, some interesting behaviors become apparent. For instance, as $n$ increases, $\ell(K_n^-)$ tends to zero, whereas $\ell(-K_n)$ tends to one (pointing to a `small' or `large' number of negative cycles, respectively). Another example is the signed graph $\G=(K_5,\sigma)$ of Figure \ref{fig1}, whose unbalance level is $\ell(\G)=0.8322>0.7912=\ell(-K_5)$ (see Table \ref{tab:ell-tokens}), indicating again a large number of negative cycles.
			Furthermore, $l(-K_5) = 4$ and 
			$l(K_5, \sigma) = 3$ and thus, there are signed graphs
			$\G=(G,\sigma)$ and $\G'=(G,\sigma')$ (with the same underlying graph $G$) 
			with $l(\G) < l(\G')$ and $\ell(\G) > \ell(\G')$.

			%\textcolor{red}{\st{ The examples of Table }\ref{tab:classesP} \st{show that there are signed graphs
					%$\G=(G,\sigma)$ and $\G'=(G,\sigma')$ (with the same underlying graph $G$) 
					%with $l(\G) < l(\G')$ and $\ell(\G) > \ell(\G')$.} \blue{(Better now with the `adjusted' definition of the unbalance level)} 
				%Furthermore, 
				%there are signed graphs
				%$\G_1=(G,\sigma_1)$ and $\G_1'=(G,\sigma_2)$ with 
				%$l(\G_1) = l(\G_1')$ and $\ell(\G) \not = \ell(\G')$.
				%\blue{(As commented, this is because the unbalance level seems to have a better discernment capacity)}\\
				%\mbox{[}QUESTION: Could it be true that
				%if $\ell(\G) = \ell(\G')$, then $l(\G) = l(\G')$?]\\
				
				The following problems on the relation between the frustration index and the unbalanced level with regard to token signed graphs are natural. 
				
				\begin{problem} 
					\label{Prob: l and ell monotone for token}
					Let $k \geq 1$ be an integer and $\G=(G,\sigma)$ and $\G'=(G,\sigma')$ be two signed graphs (with the same underlying graph $G$). 
					Prove or disprove the following statements:
					\begin{enumerate}
						\item   
						$\ell(\G) \leq \ell(F_k(\G))$.
						%\item 
						%\blue{\st{If $l(\G) \leq l(\G')$, then $\ell(\G) \leq \ell(\G').$ (iff?)}}
						\item 
						If $l(\G) \leq l(\G')$, then $l(F_k(\G)) \leq l(F_k(\G')).$
						\item 
						If $l(\G) \leq l(\G')$, then $\ell(F_k(\G)) \leq \ell(F_k(\G')).$
					\end{enumerate}
				\end{problem}
				Corollary \ref{coro:F_kbalanced} implies that all statements of Problem \ref{Prob: l and ell monotone for token} are true if $\G$ is balanced.  
				
				\begin{lemma}
					Let $\G=(G,\sigma)$ and $\G'=(G',\sigma')$ be two signed graphs (with the same underlying graph $G=G'$). If $\G$ and $\G'$ are cospectral, then they have the same unbalance level, $\ell(\G)=\ell(\G')$. 
				\end{lemma}
				
				\begin{proof}
					Let $\G=(G,\sigma)$ and $\G'=(G',\sigma')$ have (signed) adjacency matrices $\A$ and $\A'$. Let $\A^+= (\A')^+$ be the (unsigned) adjacency matrices of $G=G'$. Then, since both pair of matrices $(\A,\A')$ and $(\A^+,(\A')^+)$ have the same spectra, the result follows since $\tr \A^r=\tr (\A')^r$ and $\tr (\A^+)^r=\tr ((\A')^+)^r$ for every $r\ge 0$.
				\end{proof}
				Since, as commented in Section 2, switching isomorphic (and, hence, switching equivalent) signed graphs are cospectral, we have the following consequence.
				\begin{corollary}
					If $\G=(G,\sigma)$ and $\G'=(G',\sigma')$ are switching isomorphic, then they have the same unbalance level.
				\end{corollary}
				\begin{proof}
					Just notice that, in this case, $\A^+$ and $(\A')^+$ are similar.
				\end{proof}
				
				\begin{table}[t]
					\centering
					\begin{tabular}{|c|c|c|}
						\hline
						$n$ & $\ell(C_n^-)$ & $\ell(-C_n)$ \\[.1cm]
						\hline
						%  2 & $0$ & $0$\\
						3 & $\frac{2}{5}=0.4$ & $\frac{2}{5}=0.4$  \\[.1cm]
						4 & $\frac{2}{9}\approx 0.2222$ & $0$ \\[.1cm]
						5 & $\frac{2}{11}\approx 0.1818$& $\frac{2}{11}\approx 0.1818$\\[.1cm]
						6 & $\frac{2}{29}\approx 0.06897$ & $0$ \\[.1cm]
						7 & $\frac{2}{31}\approx 0.06452$ & $\frac{2}{31}\approx 0.06452$\\[.1cm]
						8 & $\frac{2}{99}\approx 0.02020$ & $0$ \\[.1cm]
						9 & $\frac{2}{101}\approx 0.01980$ & $\frac{2}{101}\approx 0.01980$\\[.1cm]
						10 & $\frac{2}{351}\approx 0.005698$ & $0$ \\[.1cm]
						11 & $\frac{2}{353}\approx 0.005666$ & $\frac{2}{353}\approx 0.005666$\\[.1cm]
						12 & $\frac{2}{1275}\approx 0.001569$ & $0$ \\[.1cm]
						13 & $\frac{2}{1277}\approx 0.001566$ & $\frac{2}{1277}\approx 0.001566$\\[.1cm]
						14 & $\frac{2}{4707}\approx 0.0004249$ & $0$\\[.1cm]
						15 & $\frac{2}{4709}\approx 0.0004247$ & $\frac{2}{4709}\approx 0.0004247$\\[.1cm]
						\hline
					\end{tabular}
					\caption{The unbalance levels of the cycle $C_n^-$ with only one negative edge and the all-negative cycle $-C_n$.}
					\label{tab:all(-Cn)}
				\end{table}     
				
				\begin{table}[t] 
					\centering
					\begin{tabular}{|c|c|c|c|}
						\hline
						$n$ &  $\ell(K_n^-)$ & $\ell(-K_n^+)$ & $\ell(-K_n)$\\[.1cm]
						\hline
						2 & $0$ & $0$ & 0\\
						3 & $\frac{2}{5}=0.4$ & $0$ & $\frac{2}{5}=0.4$ \\[.1cm]
						4 & $\frac{3}{7}\approx 0.4286$ & $\frac{3}{7}\approx 0.4286$ & $\frac{33}{41}\approx 0.8049$\\[.1cm]
						5 & $\frac{132}{323}\approx 0.4087$ & $\frac{260}{323}\approx 0.8050$ & $\frac{72}{91}\approx 0.7912$\\[.1cm]
						6 & $\frac{191}{623}\approx 0.3066$ & $\frac{432}{511}\approx 0.8454$ &$\frac{3685}{3907}\approx 0.9432$\\[.1cm]
						7 & $\frac{71530}{264393}\approx 0.2705$ & $\frac{234460}{26393}\approx 0.8868$ & $\frac{41130}{47989}\approx 0.8571$\\[.1cm]
						8 & $\frac{160479}{680222}\approx 0.2359$ & $\frac{89025}{99481}\approx 0.8949$ & $\frac{930769}{960801}\approx 0.9687$\\[.1cm]
						9 & $\frac{8882272}{42248291}\approx 0.2102$
						& $\frac{38378032}{42248291}\approx 0.9084$ & $\frac{15149792}{17043521}\approx 0.8889$\\[.1cm]
						10 & $\frac{30438745}{164500272}\approx 0.1850$ & $\frac{154267491}{168987244}\approx 0.9129$ & $\frac{427131081}{ 435848051}\approx 0.9800$\\[.1cm]
						11 & $\frac{419505858}{2429496991}\approx 0.1727$ & $\frac{29137900440}{31583460883}\approx 0.9226$
						& $\frac{3060912150}{3367003367}\approx 0.9091$\\[.1cm]
						12 & $\frac{39380760975}{248308506892}\approx 0.1586$ & $\frac{7055280525}{7610760577}\approx 0.9270$ &$\frac{309483909361}{313842837673}\approx 0.9861$\\[.1cm]
						13 & $\frac{14930359919508}{101789500494025}\approx 0.1467$ & 
						$\frac{18993093316692}{20357900098805}\approx 0.9330$ &$\frac{8287800739416}{8287800739416}\approx 0.9231$\\[.1cm]
						14 & $\frac{64019548410273}{469166837153174}\approx 0.1365$ & 
						$\frac{136351660121952}{145601300898805}\approx 0.9365$ & $\frac{324766589636749}{328114698808275}\approx 0.9898$\\[.1cm]
						15 & $\frac{6318610394831878}{49526346984895587}\approx 0.1276$ & $\frac{46599262432248812}{49526346984895587}\approx 0.9409$ & $\frac{10424392044215786}{ 11168991475945493}\approx 0.9333$\\[.1cm]
						\hline
					\end{tabular}
					\caption{The unbalanced levels of the complete graph $K_n^-$ with only one negative edge, the negative complete graph $-K_n^+=-K_n^-$ with only one positive edge, and the all-negative complete graph $-K_n$.}
					\label{tab:all(-Kn)}
				\end{table}          
				
				\begin{table}[t]
					\centering
					\begin{tabular}{|c|c|}
						\hline
						$n$ &  $\spec(K_n^-)$ \\[.1cm]
						\hline
						2 & $1,-1$ \\[.1cm]
						3 & $1^{[2]},-2$  \\[.1cm]
						4 & $1,\pm\sqrt{5},-1$ \\
						5 & $1,\frac{1}{2}(1\pm\sqrt{33}),-1^{[2]}$ \\[.1cm]
						6 & $1,1\pm \sqrt{12},-1^{[3]}$\\[.1cm]
						7 & $1,\frac{1}{2}(3\pm\sqrt{65}),-1^{[4]}$\\[.1cm]
						8 & $1,2\pm \sqrt{21},-1^{[5]}$\\[.1cm]
						9 & $1,\frac{1}{2}(5\pm\sqrt{105}),-1^{[6]}$\\[.1cm]
						10 & $1,3\pm \sqrt{31},-1^{[7]}$\\[.1cm]
						11 & $1,\frac{1}{2}(7\pm\sqrt{153}),-1^{[8]}$\\[.1cm]
						12 & $1,4\pm\sqrt{45},-1^{[9]}$\\[.1cm]
						13 & $1,\frac{1}{2}(9\pm\sqrt{209}),-1^{{[10]}}$\\[.1cm]
						14 & $1,5\pm \sqrt{60},-1^{{[11]}}$\\[.1cm]
						15 & $1,\frac{1}{2}(11\pm\sqrt{273}),-1^{{[12]}}$\\[.1cm]
						\hline
					\end{tabular}
					\caption{The spectrum of the complete graph $K_n^-$ with only one negative edge.}
					\label{tab:sp(Kn^-)}
				\end{table} 
				
				\begin{table}[t]
					\centering
					\begin{tabular}{|c|c|c|}
						\hline
						Graph &  Frustration index & Unbalance level \\[.1cm]
						\hline
						$+P \simeq -P_{3,3}$ & $0$ & $0$ \\[.1cm]
						$P_1 \simeq -P_{2,3}$ & $1$ &  $\frac{752}{2069}\approx 0.3635$ \\[.1cm]
						$P_{2,2}\simeq -P_{2,2}$ & $2$ &  $\frac{5536}{8569}\approx 0.6460$\\[.1cm]
						$P_{2,3}\simeq -P_{1}$ & $2$ & $\frac{6904}{10345}\approx 0.6674$\\[.1cm]
						$P_{3,2}\simeq -P_{3,2}$ & $3$ & $\frac{168}{235}\approx 0.7149$\\[.1cm]
						$P_{3,3}\simeq -P$ & $3$ & $\frac{1944}{2821}\approx 0.6891$\\[.1cm]
						\hline
					\end{tabular}
					\caption{The frustration index and the unbalance level of the six switching isomorphic classes of the Petersen graph.}
					\label{tab:classesP}
				\end{table}  
				
				\begin{table}[t]
					\centering
					\begin{tabular}{|c|c|c|}
						\hline
						Graph &  $\ell(\G)$ & $\ell(F_2(\G))$ \\
						[.1cm]
						\hline 
						$\G$ (Fig. \ref{fig1}) & $\frac{124}{149}\approx 0.8322$ & $\frac{601808}{607349}\approx 0.9909$ \\[.1cm]
						$\G$  (Fig. \ref{fig2}) & $\frac{152}{261}\approx 0.5824$ &  $\frac{249552}{321191}\approx 0.7770$ \\[.1cm]
						$\G=C_5^-$  & $\frac{2}{11}\approx 0.1818$ &  $\frac{59}{96}\approx 0.6146$ \\[.1cm]
						$\G=-K_5$   & $\frac{72}{91}\approx 0.7912$ &  $\frac{1036224}{1209157}\approx 0.8570$ \\[.1cm]
						$\G$  (Fig. \ref{fig3}) & $\frac{35}{68}\approx 0.5147$ &  $\frac{1207426567618449571}{1216092377313656965}\approx 0.9929$ \\[.1cm]
						$\G$  (Fig. \ref{fig4}) & $0$ &  $0$ \\[.1cm]
						$\G$  (Fig. \ref{fig4} with only $-\{2,3\}$) & $\frac{1}{3}\approx 0.3333$ &  $\frac{22}{39}\approx 0.5641$ \\[.1cm]
						\hline
					\end{tabular}
					\caption{The frustration indices of some signed graphs and their 2-token signed graphs.}
					\label{tab:ell-tokens}
				\end{table}  
				
				In order to compare the frustration index with the unbalance level, Table \ref{tab:classesP} shows the two parameters for the mentioned six switching isomorphic class of the Petersen graph. As done by Zaslavsky \cite{z12}, we ordered them by increasing the order of their frustration index. Moreover, as in the previous examples with cycles and complete graphs, the unbalanced level is shown to be a fine measure to distinguish between different switching isomorphic classes of signed graphs.
				
				%%%%%%%%%%%%%%%%%%%%%%%%%%%%%%%%%%%%%%%%%%%%%%%%%%%%%%%%
				\section{Sign-symmetric token graphs} 
				\label{sec:sign-symm}
				
				If a signed graph $\G=(G,\sigma)$ and its negation $-\G=(G,-\sigma)$ are switching isomorphic, then their adjacency matrices $\A$ and $-\A$ are similar. Then, the spectrum of $\G$ is symmetric with respect to the origin
				(a characterization of bipartiteness for unsigned graphs).
				Moreover, since clearly $F_k(-\G)=-F_k(\G)$, the following result holds.
				
				\begin{lemma}
					\label{lem:G-noG}
					If a signed graph $\G$ and its negation $-\G$ are switching isomorphic,  so are the  $k$-token signed graph $F_k(\G)$ and its negation $-F_k(\G)$. 
				\end{lemma}
				
				As a consequence, we have that $k$-tokens graphs preserve sign-symmetry.
				\begin{theorem}
					If $\G$ is a sign-symmetric graph, then its $k$-token signed graph $F_k(\G)$ is also sign-symmetric.
				\end{theorem}
				\begin{proof}
					From the hypothesis, $\G$ is switching isomorphic to its negation $-\G$. Then, by Theorem \ref{theo:switch-equiv} and Lemma \ref{lem:G-noG}, the $k$-token graphs $F_k(\G)$ and $-F_k(\G)$ are also switching isomorphic.
				\end{proof}
				
				An example of a sign-symmetric graph and its sign-symmetric 2-token graph is shown in Figure \ref{fig3}. As commented above, both signed graphs have spectra that are symmetric with respect to the origin. Indeed, such spectra are the following (with approximated values):
				\begin{eqnarray*} 
					\spec(\G) & = &\{0^{[2]},\pm 1, (\pm\sqrt{5})^{[2]}\};\\
					\spec(F_2(\G)) & = & \{0^{[6]},\pm 0.8140,\pm 1,\pm 1.236,\pm 1.448, \pm 1.536, \pm 2.236,\\
					& & \pm 2.628, \pm 2.888, \pm 3, \pm 3.236, \pm4.318\}.
				\end{eqnarray*}

				\begin{figure}[t]
					\begin{center}
						\vskip .5cm
						\includegraphics[width=10cm]{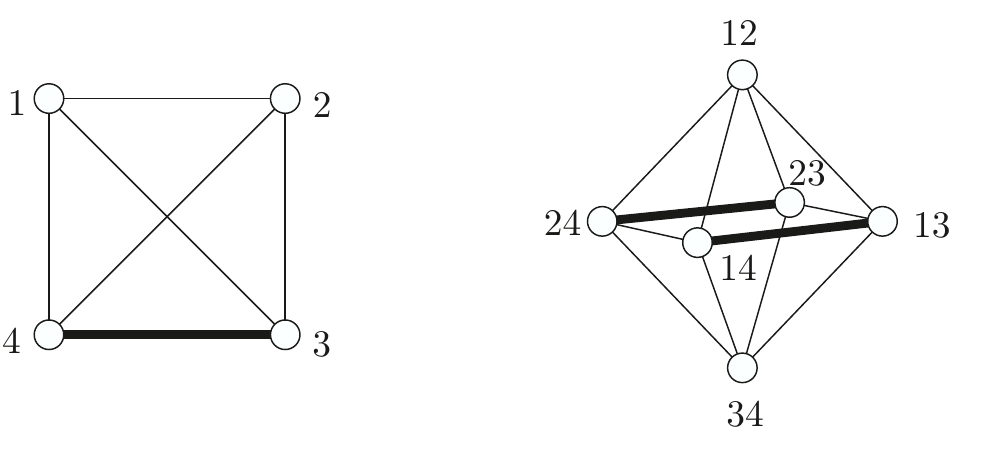}
						\caption{The sign-symmetric signed graph $K_4^-$ and its sign-symmetric 2-token signed graph.}
						%, \textcolor{blue}{which is balanced if $V_1=\{12,23,34,45,51\}$ and $V_2=\{13,24,35,42,52\}$}. 
						\label{fig:K4}
					\end{center}
				\end{figure}
				
				Another example is the complete graph $\G=K_4^-$ (see Figure \ref{fig:K4}) with only one negative edge satisfying (see Table \ref{tab:all(-Kn)} where $\ell(K_4^-)=\ell(-K_4^+$)). Its spectrum and that of its 2-token graph are
				$$
				\spec(\G)=\{\pm 1, \pm\sqrt{5}\}\quad \mbox{and}\quad \spec(F_2(\G))=\{0^{[2]},\pm 2,\pm 2\sqrt{2}\}.
				$$
				
				\begin{figure}[t]
					\begin{center}
						\includegraphics[width=10cm]{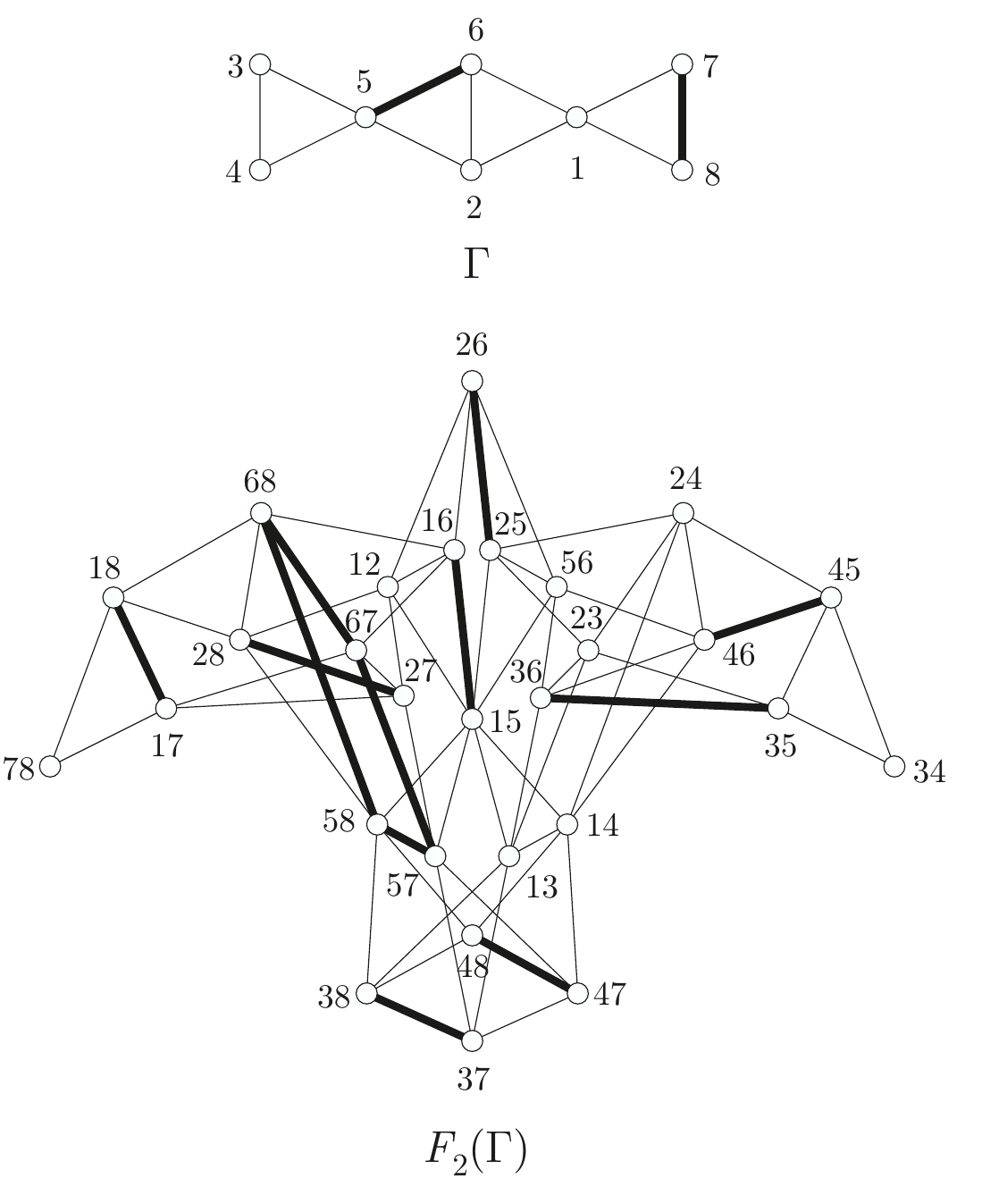}
						%\vskip -1.5cm
						\caption{A sign-symmetric graph and its sign-symmetric 2-token graph (the `bird graph').}
						\label{fig3}
					\end{center}
				\end{figure}
				% \begin{figure}[t]
					% 	\begin{center}
						% 		\includegraphics[width=9cm]{signed-bird-graph}
						% %\vskip -1.5cm
						% 	\caption{The signed bird graph ($F_2(\G)$ in Figure \ref{fig3}).}
						% \label{fig4}
						% 	\end{center}
					% \end{figure}
				
				%%%%%%%%%%%%%%%%%%%%%%%%%%%%%%%%%%%%%%%%%%%%%%
				\section{The Laplacian spectrum}
				\label{sec:Lapl}
				
				Given some integers $n$ and $k$, with $k\in [1,n-1]$, let $\G=(G,\sigma)$ a signed graph on $n$ vertices. Then, the {\em signed $(n,k)$-binomial matrix} $\B(\G)$ is an ${n\choose k}\times n$ matrix whose rows are indexed by the $k$ subsets $A_1,\ldots, A_{{n\choose k}}$ of $[n]$, and its entries are
				$$
				(\B(\G))_{ij}=\left\{
				\begin{array}{cc}
					\pm 1 & \mbox{if } j\in A_i,\\
					0 & \mbox{otherwise},
				\end{array}
				\right.
				$$
				where the plus-minus sign of $(\B(\G))_{ij}$ depends on $\G$, as shown in the following result.
				
				\begin{proposition}
					\label{BL1=L2B}
					Let $\G=(G,\sigma)$ be a balanced signed graph on $n$ vertices, and $F_k(\G)$ its signed $k$-token graph for some $k\in [1,n-1]$. Let $\L_1$ and $\L_k$ be their corresponding Laplacian matrices. Then, there exists a signed binomial matrix $\B$ such that
					\begin{equation}
						\B\L_1=\L_k\B.
						\label{eq:BL1=LkB}
					\end{equation}
				\end{proposition}
				
				\begin{proof}
					Let $\L_1$ and $\L_k$ be the Laplacian matrices of the underlying graphs $G$ and $F_k(G)$, respectively.
					Since $\G$ and $F_k(\G)$ are balanced, they are switching equivalent to the unsigned graphs  $G$ and $\L_k(G)$, with respective Laplacian matrices $\L_1^+$ and $\L_k^+$. Let $U$ be the switching set that leads from $\G$ to $G$. Then, the set $U_k$, constructed from $U$ as in the proof of Theorem \ref{theo:switch-equiv},
					corresponds to the switching equivalence between $F_k(\G)$ and $F_k(G)$. Let $\S_1$ and $\S_k$ be the diagonal $(+1,-1)$-matrices representing $U$ and $U_k$, so that 
					$$
					\L_1^+=\S_1\L_1\S_1\qquad\mbox{and}\qquad \L_k^+=\S_k\L_k\S_k.
					$$
					By Theorem \ref{theo:BL1=LB2}, we have
					$$
					\B^+\L_1^+=\L_k^+\B^+,
					$$
					where $\B^+$ is the standard (unsigned) $(n;1,k)$-binomial matrix.
					Thus,
					$$
					\B^+\S_1\L_1\S_1=\S_k\L_k\S_k\B^+
					$$
					or
					$$
					\S_k^{-1}\B^+\S_1\L_1 =\L_k\S_k\B^+\S_1^{-1}.
					$$
					However, since $\S_1$ and $\S_k$ are diagonal $(+1,-1)$-matrices, we have that they coincide with their inverses.
					Consequently, the signed binomial matrix $\B=\S_k\B^+\S_1$
					satisfies  \eqref{eq:BL1=LkB}.
					%Eq. \eqref{eq:BLk1=Lk2B} in Theorem \ref{BLk1=Lk2B}.  
				\end{proof}
				
				As a consequence of Proposition \ref{BL1=L2B}, we have the following result.
				\begin{theorem}
					Let $\G=(G,\sigma)$ be a balanced signed graph on $n$ vertices, and $F_k(\G)$ its signed $k$-token graph for some $k\in [1,n-1]$. Then, the Laplacian spectrum of $\G$ is contained in the Laplacian spectrum of $F_k(G)$.
				\end{theorem}
				
				\begin{proof}
					Let $\vv$ be an eigenvector of $\L_1$ with eigenvalue $\lambda$. Then, we claim that $\B\vv$ is an eigenvector of $\L_k$ with the same eigenvalues $\lambda$. Indeed, from $\L_1\vv=\lambda\vv$ and Proposition \ref{BL1=L2B}, we have that
					$$
					\L_k\B\vv=\B\L_1\vv=\lambda\B\vv.
					$$
					Moreover, independent eigenvectors of $\L_1$ give rise to independent eigenvectors of $\L_k$. The reason is that, since  $\rank(\B^+)=n$ (see de Caen \cite{dc01}), we also have $\rank(\B)=n$ and, hence, $\Ker(\B)=\{\vec0\}$.
				\end{proof}
				
\begin{figure}[t]
					\begin{center}
						\includegraphics[width=8cm]{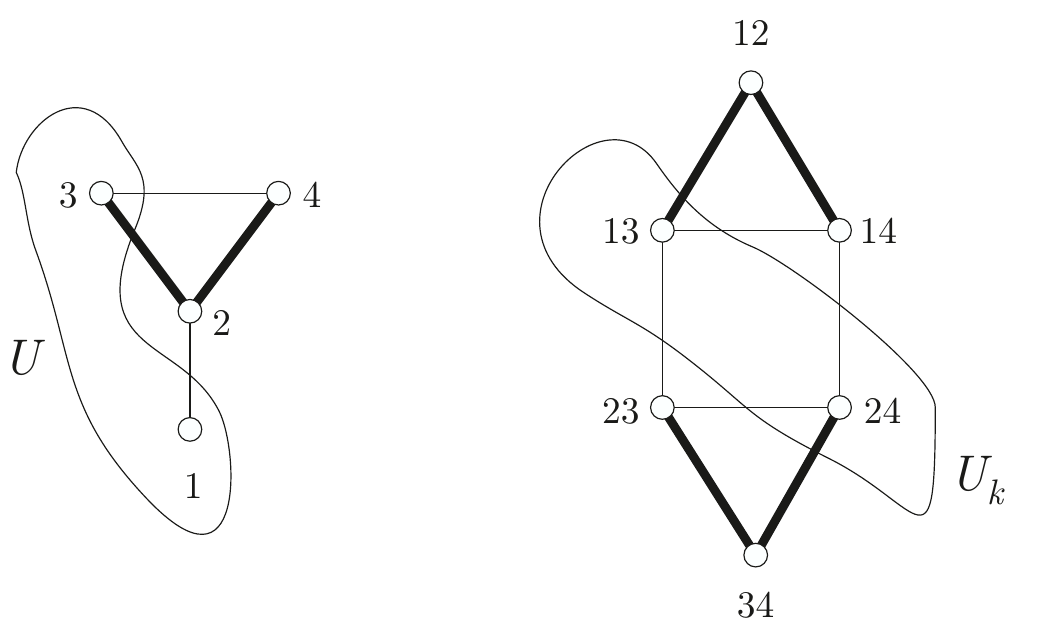}
						%\vskip -1.5cm
						\caption{A balanced signed  graph $\G$ and its balanced  2-token signed graph $F_2(\G)$.}
						\label{fig4}
					\end{center}
				\end{figure}

				For instance, if $\G$ and $F_2(\G)$ are the balanced signed graphs shown in Figure \ref{fig4}, their respective signed Laplacian matrices are
				$$
				\L_1=\left( 
				\begin{array}{cccc}
					1& -1& 0& 0\\
					-1& 3& 1& 1\\
					0& 1& 2& -1\\
					0& 1& -1& 2
				\end{array}
				\right)
				\quad \mbox{and}\quad
				\L_2=\left( 
				\begin{array}{cccccc}
					2& 1&  1&  0&  0&  0\\
					1&  3&  -1&  -1&  0&  0\\
					1&  -1&  3&  0&  -1&  0\\
					0&  -1&  0&  3&  -1&  1\\
					0&  0&  -1&  -1&  3&  1\\
					0&  0&  0&  1&  1&  2
				\end{array}
				\right).
				$$
				They are switching isomorphic to all positive (underlying) graphs, with corresponding matrices $\S_1=\diag(+1,-1,+1,-1)$ and 
				$\S_2=\diag(-1,+1,-1,-1,$ $+1,-1)$. Thus,  the signed binomial matrix satisfying Proposition \ref{BL1=L2B} turns out to be
				$$
				\B=\S_2\B^+\S_1=\left( 
				\begin{array}{cccc}
					1& 1& 0& 0\\
					-1& 0& 1& 0\\
					-1& 0& 0& 1\\
					0& -1& 1& 0\\
					0& -1& 0& 1\\
					0& 0& -1& -1
				\end{array}
				\right).
				$$
				Moreover, the spectra  of $\L_1$  and $\L_2$ are
				$$
				\spec (\L_1)=\{0^{[2]},1^{[2]},3^{[2]},4^{[2]}\} \subset 
				\{0^{[2]},1^{[2]},3^{[4]},4^{[2]},5^{[2]}\}=\spec (\L_2). 
				$$
				In fact, using Theorem \ref{theo:BL1=LB2}, we have the following generalization of Proposition \ref{BL1=L2B} (the proof is basically the same as that of this result).
				\begin{theorem}
					\label{BLk1=Lk2B}
					Let $\G=(G,\sigma)$ be a balanced signed graph on $n$ vertices. Given integers $1\le k_1<k_2<n$, let  $F_{k_1}(\G)$ and  $F_{k_2}(\G)$ be their $k_1$-token and $k_2$-token signed graphs, with corresponding Laplacian matrices $\L_{k_1}$ and $\L_{k_2}$. Then, there exists a signed $(n;k_2,k_1)$-binomial matrix $\B$ such that
					\begin{equation}
						\B\L_{k_1}=\L_{k_2}\B.
						\label{eq:BLk1=Lk2B}
					\end{equation}
				\end{theorem}
				
				%%%%%%%%%%%%%%%%%%%%%%%%%%%%%%%%%%%%%%%%%%%%%%%%%%%%%%%%%%%%% 
				\section{The complement signed graph and its $k$-token graph}
				\label{sec:compl}
				
				Let $G$ be an unsigned graph with $n$ vertices, adjacency matrix $\A$, and Laplacian matrix $\L=\D-\A$. Let $\overline{G}$ be the complement of $G$, with adjacency and Laplacian matrices $\overline{\A}=\J-\I-\A$ and $\overline{\L}=\overline{\D}-\overline{\A}=n\I-\D-\J+\A=n\I-\J-\L$.
				
				Let $\G_{\sigma}=(G,\sigma)$ be a balanced signed graph with adjacency matrix $\A_{\sigma}$ and Laplacian matrix $\L_{\sigma}$. Since $\G$ is switching equivalent to the underlying graph $G$, there exists a diagonal $(+1,-1)$-matrix $\S$ such that $\A_{\sigma}=\S\A\S$ and $\L_{\sigma}=\S\L\S$.  This allows us to define the complement of a balanced signed graph.
				\begin{definition}
					\label{def:comp}
					Let $\G_{\sigma}=(G,\sigma)$ be a balanced signed with matrices $\A$ and $\L$, and diagonal matrix $\S$ defined above.
					The {\em signed complement} of  $\G_{\sigma}$, denoted by $\overline{\G}_{\sigma}$, is the signed graph with adjacency matrix $\overline{\A}_{\sigma}=\S\overline{\A}\S$ and Laplacian matrix $\overline{\L}_{\sigma}=\S\overline{\L}\S.$
				\end{definition}
				
				\begin{proposition}
					\label{propo:complement}
					Let $\G_{\sigma}=(G,\sigma)$ be a balanced signed graph on $n$ vertices, and $\overline{\G}_{\sigma}=(\overline{G},\overline{\sigma})$ its signed complement graph. Then, the following statements hold.
					\begin{itemize}
						\item[$(i)$]
						The signed graph $\overline{\G}_{\sigma}$ is also balanced.
						\item[$(ii)$]  
						The matrices $\M_A=\A_{\sigma}+\overline{\A}_{\sigma}$
						and $\M_L=\L_{\sigma}+\overline{\L}_{\sigma}$ are, respectively, the adjacency matrix and the Laplacian matrix of a balanced signed complete graph ${\cal K}_n=(K_n,\tau)$.
						\item[$(iii)$] 
						The Laplacian matrices $\L_{\sigma}$ and $\overline{\L}_{\sigma}$
						commute:
						$$
						\L_{\sigma}\overline{\L}_{\sigma}=\overline{\L}_{\sigma}\L_{\sigma}.
						$$
						\item[$(iv)$] 
						Let $\lambda_0(=0)\le \lambda_1\le \cdots\le \lambda_{n-1}$
						and $\overline{\lambda}_0(=0)\le \overline{\lambda}_1\le \cdots\le \overline{\lambda}_{n-1}$ be the eigenvalues of 
						$\L_{\sigma}$ and $\overline{\L}_{\sigma}$, respectively. Then,
						$$
						\lambda_i+\overline{\lambda}_{n-i}=n,\quad \mbox{ for \ $i=1,\ldots,n-1$.}
						$$
					\end{itemize}
				\end{proposition}
				
				\begin{proof}
					$(i)$ By Definition \ref{def:comp}, $\overline{\G}_{\sigma}$ is switching equivalent to the (balanced) complement graph $\overline{G}$.
					
					\item[$(ii)$]
					The matrices $\M_A$ and $\M_L$ are the signed matrices of $\A+\overline{\A}$
					and $\L+\overline{\L}$, that correspond to the adjacency and Laplacian matrices of $K_n$. Again, from Definition \ref{def:comp}, ${\cal K}_n$ is switching equivalent to the (balanced) graph $K_n$.
					
					$(iii)$ Let $\S$ be the switching matrix in Definition \ref{def:comp}. Since $\L+\overline{\L}=n\I-\J$ and $\L$ and $\overline{\L}$ commute with  $\I$ and $\J$, we have that $\L\overline{\L}=\overline{\L}\L$. Since  $\S^2=\I$, so that
					$$
					\S\L\S^2\overline{\L}\S=\S\overline{\L}\S^2\L\S
					$$
					gives the result.
					
					$(iv)$ This is a consequence of $(ii)$ and $(iii)$. First, notice that, since the matrix $\M_L=\L_{\sigma}+\overline{\L}_{\sigma}$ in $(i)$ is switching equivalent to the Laplacian matrix of $K_n$, its eigenvalues are $0$ (with multiplicity 1) and $n$ (with multiplicity $n-1$). 
					Moreover, by a theorem of Frobenius, since the matrices $\L_{\sigma}$ and $\overline{\L}_{\sigma}$ commute, their eigenvalues
					can be matched up as
					$\lambda_i\leftrightarrow \overline{\lambda}_i$ in such a way that the eigenvalues of any polynomial $p(\L_{\sigma},\overline{\L}_{\sigma})$ in the two matrices is the multiset of the values $p(\lambda_i,\overline{\lambda}_i)$.
				\end{proof}
				
				Notice that neither $(ii)$ nor $(iii)$ hold for the adjacency matrices $\A_{\sigma}$ and $\overline{\A}_{\sigma}$.
				
				Finally, Corollary \ref{coro:F_kbalanced} and Proposition \ref{propo:complement} lead us to the following result.
				
				\begin{theorem}
					\begin{itemize}
						\item[$(i)$]
						Let $\G$ be a balanced signed graph and let $\overline{\G}$ be its signed complement, both on $n$ vertices. The Laplacian (signed) matrices $\L_k=\L(F_k(\G))$ and $\overline{\L}_k=\L(F_k(\overline{G}))$ of their $k$-token signed graphs commute:
						$$
						\L_k\overline{\L}_k=\overline{\L}_k \L_k.
						$$
						\item[$(ii)$]
						Let $\G$ and $\overline{\G}$ be as above. 
						Then, every eigenvalue $\lambda_J$ of the Johnson graph $J(n,k)=F_k(K_n)$ is the sum of one eigenvalues $\lambda_{F_k(\G)}$ of $F_k(\G)$ and one eigenvalue $\lambda_{F_k(\overline{\G})}$ of $F_k(\overline{\G})$, where each $\lambda_{F_k(\G)}$ and each $\lambda_{F_k(\overline{\G})}$ is used once:
						$$
						\lambda_{F_k(\G)}+\lambda_{F_k(\overline{\G})}=\lambda_J.
						$$
					\end{itemize}
				\end{theorem}

				% %%%%%%%%%%%%%%%%%%%%%%%%%%%%%%%%%%%%%%%%%%%%%
				% \section{Problems}
				
				% \begin{problem}
					% Define the (signed) partial $k$-token of a signed graph in order to construct frustration-critical signed graphs.
					% \end{problem}

				% \begin{problem}
					% Study tokens of signed bipartite graphs together with their signless Laplacian.
					% (The signless Laplacian of a bipartite graph has the same spectrum as that of its Laplacian matrix.)
					% In particular, study $k$-tokens of $K_{n,n}$ as a variation of Johnson graphs.
					
					% \begin{problem}
						% Find a spectral measure of closeness to a balanced signed graph.
						% \end{problem}
					% \end{problem}
				
				%%%%%%%%%%%%%%%%%%%%%%%%%%%%%%%%%%%%%%%%%%%%%%%%%%%%%%%%%%%
				\section{Declarations}
				\subsection{Data availability}
				There is no data associated with this paper. \\
				\subsection{Funding and/or Conflicts of interests/Competing interests}
				The authors do not have any conflict of interest or competing interests.

				%%%%%%%%%%%%%%%%%%%%%%%%%%%%%%%%%%%%%%%%%%%%%%%%%%%%%%%%%%%%%

			\end{document}